\documentclass[reqno]{amsart}

\usepackage {amsmath, amssymb}
\usepackage {fancybox}
\usepackage [top=1.5in, right=1in, bottom=1.5in, left=1in] {geometry}

\usepackage {color}

\title {Log-Canonical Coordinates for Poisson Brackets and Rational Changes of Coordinates}

\author{John Machacek}
\address{Department of Mathematics, Michigan State University, USA}
\email{machace5@math.msu.edu}

\author{Nicholas Ovenhouse}
\address{Department of Mathematics, Michigan State University, USA}
\email{ovenhou3@math.msu.edu}

\keywords{Poisson algebras, Poisson varieties, log-canonical coordinates, cluster algebras}
\subjclass[2010]{Primary 17B63; Secondary 	13F60, 53D17}

\newcommand{\Z}{\mathbb{Z}}
\newcommand{\N}{\mathbb{N}}
\newcommand{\supp}{\mathrm{supp}}
\newcommand{\ad}{\mathop{\mathrm{ad}}}

\newcommand{\Lau}{\mathcal{L}}
\newcommand{\Rat}{\mathcal{R}_{\Omega}}
\newcommand{\K}{\mathbb{K}}
\newcommand{\I}{\mathcal{I}}
\newcommand{\J}{\mathcal{J}}

\newcommand{\x}{\mathbf{x}}

\theoremstyle {definition}
\newtheorem*  {ex} {Example}
\newtheorem*  {rem} {Remark}

\theoremstyle {plain}
\newtheorem  {thm} {Theorem}
\newtheorem* {thm*} {Theorem}
\newtheorem {prp}[thm] {Proposition}
\newtheorem  {lem}[thm] {Lemma}
\newtheorem  {cor}[thm] {Corollary}
\newtheorem{conj}[thm]{Conjecture}
\newtheorem{question}[thm]{Question}

\begin {document}

\begin{abstract}
Goodearl and Launois have shown in \cite{GL11} that for a log-canonical Poisson bracket on affine space there is no rational change of coordinates for which the Poisson bracket is constant.
Our main result is a proof of a conjecture of Michael Shapiro which states that if affine space is given a log-canonical Poisson bracket,
then there does not exist any rational change of coordinates for which the Poisson bracket is linear.
Hence, log-canonical coordinates can be thought of as the simplest possible algebraic coordinates for affine space with a log-canonical coordinate system.
In proving this conjecture we find certain invariants of log-canonical Poisson brackets on affine space which linear Poisson brackets do not have.
\end{abstract}
\maketitle

\tableofcontents

\section{Introduction}

Cluster algebras were originally defined by Fomin and Zelevinsky~\cite{FZ02} to study total positivity and dual canonical bases in semisimple groups.
Since their inception, connections between cluster algebras and many areas of algebra and geometry have been found.
One such connection is with Poisson geometry.
Gekhtman, Shapiro, and Vainshtein~\cite{GSV10} have studied Poisson structures compatible with cluster algebras.
In this compatibility, the cluster variables give log-canonical coordinates for the Poisson bracket, while the mutations
give birational transformations preserving the log-canonicity.
We will study log-canonical Poisson brackets under rational changes of coordinates.
Our main result is Theorem~\ref{thm:main} where we show that log-canonical coordinates are analogous to Darboux coordinates for rational algebraic functions
in the sense that the Poisson bracket takes the simplest form in these coordinates.
  
\subsection{Poisson algebras and Poisson geometry}
Let $P$ be an associative algebra.
A \emph{Poisson bracket} on $P$ is a skew-symmetric bilinear map $\{\cdot ,\cdot \}: P \times P \to P$ such that for any $a,b,c \in P$ both the \emph{Leibnitz identity}
$$\{ab,c\} = a\{b,c\} + \{a,c\}b$$
and the \emph{Jacobi identity}
$$\{a,\{b,c\}\} + \{b,\{c,a\}\} + \{c,\{a,b\}\} = 0$$
hold.
A \emph{Poisson algebra} is pair $(P, \{\cdot, \cdot\})$ where $P$ is an associative algebra and $\{\cdot, \cdot\}$ is a Poisson bracket.

Notice that $\{\cdot,\cdot\}$ makes $P$ a Lie algebra. So, we get the adjoint representation of $P$ on itself
sending $a \in P$ to $\ad_a \in \mathrm{End}(P)$, where $\ad_a(b) = \{a,b\}$. Note that the Jacobi identity implies that $\ad_a$ is a Lie
algebra derivation. Also observe that $\ad_a$ is a derivation of the associative algebra $P$ by the Leibniz identity. If $a \in P^*$
is a unit, then the Leibniz identity implies that $\ad_{a^{-1}} = -a^{-2} \ad_a$.
In particular, this implies that if $\{a,b\} = 0$ for some $a \in P^*$ and $b \in P$, then $\{a^{-1}, b\} = 0$.

Let $M$ be a smooth manifold, and let $C^{\infty}(M)$ denote its algebra of smooth functions.
A \emph{Poisson structure} on $M$ is a bracket $\{\cdot, \cdot\}:C^{\infty}(M) \times C^{\infty}(M) \to C^{\infty}(M)$  such that $(C^{\infty}(M), \{\cdot, \cdot\})$ is a Poisson algebra.
In this case we call  $(M, \{\cdot, \cdot\})$ a  \emph{Poisson manifold}.
For local coordinates $(x_1, \dots, x_n)$ and $f,g \in C^{\infty}(M)$ the Poisson bracket is given by
\begin {equation} 
    \{f,g\} = \sum_{i,j = 1}^n \frac{\partial f}{\partial x_i} \frac{\partial g}{\partial x_j} \{x_i, x_j\}.
    \label{eq:bracket-formula}
\end {equation}
and so the bracket is completely determined by the $\binom{n}{2}$ ``structure functions'' $\{x_i,x_j\}$, for $i<j$.
Following \cite{GSV10}, a system of coordinates $(x_1, \dots, x_n)$ is called \textit{log-canonical} with respect to a Poisson bracket $\{\cdot, \cdot\}$ if 
there is a matrix of scalars $\Omega = (\omega_{ij})$ (necessarily skew-symmetric) such that the structure functions are
given by $\{x_i, x_j\} = \omega_{ij} x_i x_j$. We note here that this Poisson structure goes by many names in the literature.
For example, it is called a ``\textit{diagonal Poisson structure}'' in \cite{VH2}, ``\textit{Poisson $n$-space}'' in \cite{OH06},
and a ``\textit{semi-classical limit of quantum affine space}'' in \cite{GOODEARL}.

In general, the local structure of Poisson manifolds is described by the following theorem of Weinstein.

\begin {thm*}[\cite{W}] Let $M$ be a Poisson manifold, and $p \in M$. Then there exists a neighborhood $U$ containing $p$
with coordinates $(x_1,y_1,\dots,x_r,y_r,z_1,\dots,z_s)$, such that the bracket takes the form
\[ \{x_i,x_j\} = \{y_i,y_j\} = \{x_i,z_j\} = \{y_i,z_j\} = 0 \]
\[ \{x_i,y_j\} = \delta_{ij} \]
\[ \{z_i,z_j\} = \varphi_{ij} \]
where $\varphi_{ij} \in C^\infty(U)$ depend only on $z_1,\dots,z_s$, and $\varphi_{ij}(p) = 0$.
\end {thm*}

\bigskip

\begin {ex}
    If $(M^{2n}, \omega)$ is a symplectic manifold, then there is a standard Poisson structure induced by $\omega$.
    In this special case, Weinstein's theorem is the classical Darboux theorem which says that locally $\omega$ has the form
    \[ \omega = \sum_{i=1}^n dx_i \wedge dy_i \]
    The local coordinates $(x_1,y_1,\dots,x_n,y_n)$ are commonly called ``canonical coordinates'' or ``Darboux coordinates.''
\end {ex}

Note that on a smooth Poisson manifold with a log-canonical system of coordinates $(x_1,\dots,x_n)$ the system
of coordinates $(y_1,\dots,y_n) = (\log x_1, \dots, \log x_n)$, defined on the open set where all $x_i$ are positive,
are similar to a system of canonical coordinates in the sense that the structure functions
$$\{y_i,y_j\} = \{\log x_i, \log x_j\} = \omega_{ij}$$ are all constants. This is indeed the intuition behind the terminology ``log-canonical.''
From Theorem~\ref{thm:constant} it will follow that there does not exist any rational change of coordinates 
on any Zariski open subset such that the structure functions are constant in the new coordinates.

Similarly, let $M$ be an algebraic variety and $\mathcal{O}(M)$ its algebra of regular functions. If there is a bracket making
$\mathcal{O}(M)$ into a Poisson algebra, then we call $(M, \{\cdot, \cdot\})$ a \emph{Poisson variety}. 
Suppose there is a system of coordinates $(x_1,\dots,x_n)$ on some Zariski open subset of a Poisson variety $M$,
then the bracket is given by Equation~(\ref{eq:bracket-formula}) just as in the smooth case (see for example \cite{VH2} for details).
We wish to investigate whether such a ``simplification'' of the structure functions is possible
(analogous to the simplification in the Darboux/Weinstein Theorem, in the sense that all structure functions become lower degree polynomials), 
allowing only birational change-of-coordinates.
It is suggested/conjectured in \cite{VH} that there are not canonical coordinates in general for an arbitrary Poisson variety, 
but that no specific counterexample has been demonstrated.
In \cite{GL11}, it was shown that affine space with a log-canonical bracket is such a counterexample.
We wish to demonstrate that this same example has the additional property that no rational change of coordinates
can make the structure functions linear.
The following example is given in~\cite{VH} and demonstrates some of the nuances of the problem of finding canonical coordinates on an open set of a Poisson variety.

\begin{ex}[\cite{VH}]
Consider affine space $\mathbb{C}^2$ with coordinates $(x,y)$ and Poisson bracket given by $\{x,y\} = x$.
Viewing $\mathbb{C}^2$ as a smooth manifold, 
there is a system of canonical local coordinates $(\log x, y)$ that is \emph{not} algebraic.
However, there is also $\left(\frac{1}{x}, -xy\right)$ which is a system of canonical local coordinates that \emph{is} algebraic.
That is, a system of canonical coordinates consisting of rational functions in $x$ and $y$ defined on the Zariski-open subset $\{(x,y) : x \neq 0\}$ of the variety $\mathbb{C}^2$.
The example illustrates that there do exist Poisson varieties which admit a rational coordinate change on an open subset which
make the structure functions constant.
\end{ex}

\begin {ex}
    More generally, consider $\mathbb{C}^2$ with coordinates $(x,y)$ and Poisson bracket given by $\{x,y\} = x^a y^b$ for $(a,b) \in \N \times \N$.
    The case $(a,b)=(1,1)$ gives a system of log-cononical coordinates.
    In all other instances, we can find a system of canonical coordinates as follows:
    \begin {itemize}
        \item If $a \neq 1$ and $b \neq 1$, then $\{x^{-(a-1)}, y^{-(b-1)}\} = (a-1)(b-1)$ is a nonzero constant.
        \item If $a = 1$ and $b \neq 1$, then $\{x^{-1}, xy^{-(b-1)}\} = (b-1)$ is a nonzero constant. The case $a \neq 1$ and $b = 1$ is similar using the fact that the bracket is antisymmetric.

    \end {itemize}
    Note that the previous example is the special case when $(a,b) = (1,0)$.
    Although the specific example $(a,b) = (1,0)$ does give a birational change of coordinates,
    this is not in general true for this family of examples. For instance, when either $a$ or $b$ is
    greater than $2$, the inverse of the coordinate change is not a rational function.

    Thus for $(a,b) \neq (1,1)$ we can always find a pair of algebraically independent rational functions in two variables 
    such that the bracket between these two functions is a nonzero constant.
    It is still unclear whether this example can be generalized to dimensions higher than 2.
    It will follow from Theorem~\ref{thm:constant} that $(a,b) = (1,1)$ is the unique exception to the existence of two rational functions with nonzero constant bracket between them.
    This begs the following interesting, and more general, question. 
\end {ex}

\begin{question}
    Given a Poisson bracket whose structure functions are all (homogeneous) polynomials of a given degree,
    when is it possible to find a birational change of coordinates making the structure functions (homogeneous) polynomials of a smaller degree?
\label{q:degree}
\end{question}

The aim of the last section is to give a more technical formulation of this question.

\bigskip

Before discussing the main results, we first review some motivations for the study of log-canonical brackets.
Log-canonical coordinates arise naturally in the standard Poisson-Lie structure on a simple Lie group and in the theory of cluster algebras.
We will briefly review these connections in the following two sections.

\subsection{Connections to Cluster Algebras}
Cluster algebras were defined by Fomin and Zelevinsky in~\cite{FZ02}.
Informally, a \emph{cluster algebra} is a subalgebra of an ambient field of rational functions in $n$ variables which is defined by some combinatorial data.
The combinatorial data is known as a \emph{seed}, and it consists of a set of \emph{cluster variables} and an \emph{exchange matrix}.
The set of cluster variables is called a \emph{cluster}.
The cluster variables are some distinguished generators of the cluster algebra, while the exchange matrix gives rules for producing more generators 
of the cluster algebra through a process called \emph{mutation}.
An initial seed will be of the form $(\{x_1, \cdots, x_n\}, B)$ where $\{x_1, \cdots, x_n\}$ are the cluster variables and $B$ is the exchange matrix.
Some of the variables ($x_1, \cdots, x_m$ for some $m \leq n$) are ``mutable'', while the remaining variables are ``frozen'' (they do not
change under the mutation process).
For any $i \leq m$, the ``mutation in direction $i$'' is given by
$$(\{x_1, \cdots, x_n\}, B) \mapsto (\{x_1', \cdots, x_n'\}, B')$$ 
where $x'_j = x_j$ when $i \neq j$, and $x_i'$ is a subtraction-free rational expression in the $x_1, \cdots, x_n$.
Specifically, the expression for $x'_i$ will be of the form
$$x'_i = \frac{\alpha m_1 + \beta m_2}{x_i}$$
where $m_1$ and $m_2$ are monomials in the cluster variables $x_1, \cdots, x_n$ whose exponents come from the exchange matrix $B$, 
and $\alpha$ and $\beta$ belong to a semifield of coefficients.
We omit further details of the definition of cluster algebras (including the rule for determining the new exchange matrix $B'$) as they will not be needed here.

The process of mutation outlined above can be iterated.
Any sequence of mutations gives a new seed, consisting of cluster variables and some exchange matrix.
The cluster algebra is then the subalgebra of rational functions in the variables $x_1, \cdots, x_n$ generated by all possible cluster variables that can be obtained by mutation.
As defined in~\cite{GSV10} a Poisson bracket on the ambient field of rational functions is called ``compatible'' with the cluster algebra 
if each cluster forms a log-canonical coordinate system.
That is, for any cluster $(\{y_1, \cdots, y_n\}, C)$ the Poisson bracket $\{\cdot, \cdot\}$ must satisfy
$$\{y_i, y_j\} = \omega_{ij} y_i y_j$$
for some skew-symmetric collection of scalars $\omega_{ij}$.

Given a cluster algebra $A$, the \emph{cluster manifold} $X(A)$ is defined in~\cite{GSV10} to be a certain nonsingular part of $\mathrm{Spec}(A)$.
In this way, a Poisson bracket on the field of rational functions compatible with the cluster algebra $A$ makes $X(A)$ a Poisson variety, 
and the cluster variables are log-canonical coordinates on this Poisson variety. Our work here justifies that these are the ``nicest''
coordinates since the bracket relations cannot be brought to a ``simpler'' form. In this sense, log-canonical coordinates are analogous
to canonical/Darboux coordinates.

\subsection{Example: Poisson-Lie Groups}
A map $\varphi \colon M \to N$ between two Poisson manifolds (or Poisson varieties) is called  a \emph{Poisson map}
if the pullback map $\varphi^*$ is a homomorphism of Poisson algebras. A Lie group $G$ is called a \emph{Poisson-Lie group} if the multiplication
map $G \times G \to G$ is a Poisson map. For further details, see \cite{CP94}.

\begin{ex}
    Consider the special linear group $\mathrm{SL}_n$, with coordinates (matrix entries) $x_{ij}$. 
    The ``standard'' Poisson-Lie structure on $\mathrm{SL}_n$ is the quadratic bracket given by
    \[ \{x_{ij}, x_{k \ell} \} = c^{ij}_{k\ell} x_{i\ell}x_{kj} \]
    where the coefficients are given by
    \[ c^{ij}_{k\ell} = \frac{1}{2} \left( \mathrm{sign}(k-i) + \mathrm{sign}(\ell - j) \right)
       = \begin{cases}
            1 & \textrm{if } k>i, \, \ell>j \\
            0 & \textrm{if } k>i, \, \ell<j \\
            \frac{1}{2} & \textrm{if } k>i,j=\ell \textrm{ or } k=i,\ell>j
         \end{cases}
    \]
    For instance, when $n=2$ we have
    \[ \mathrm{SL}_2 = \left\{ \begin{pmatrix} a&b \\ c&d \end{pmatrix} : ad - bc = 1\right\} \]
    with bracket relations:
    \begin {center}
    \begin {tabular}{cc}
        $\{a,b\} = \frac{1}{2} ab$ & $\{c,d\} = \frac{1}{2} cd$ \\ [2ex]
        $\{a,c\} = \frac{1}{2} ac$ & $\{b,d\} = \frac{1}{2} bd$ \\ [2ex]
        $\{a,d\} = bc$ & $\{b,c\} = 0$
    \end {tabular}
    \end {center}
    If we consider the Borel subgroup of upper triangular matrices in $\mathrm{SL}_2$ 
    \[ B = \left\{ \begin{pmatrix} \alpha & \beta \\ 0 & \alpha^{-1}\end{pmatrix} \right\}, \]
    then the bracket is given by $\{\alpha,\beta\} = \frac{1}{2} \alpha \beta$. In particular,
    the standard Poisson-Lie structure gives a log-canonical bracket on the Borel subgroup.
\end{ex}

\subsection{The Main Results}

We will be interested in Poisson algebras of rational functions. 
Let $\K$ be a field and $\Omega = (\omega_{ij})$ a skew-symmetric matrix. 
Consider the algebra $\Rat = \K(x_1, \dots, x_n)$ of rational functions in $n$ variables with a Poisson bracket in which
the functions $x_1,\dots,x_n$ form a system of log-canonical coordinates:
\[ \{x_i, x_j\} = \omega_{ij} x_i x_j \]
Here we wish to show that the bracket $\{\cdot, \cdot\}$ has the simplest expression in the coordinates $x_1, \cdots, x_n$.
In particular, we want to show that no rational change of coordinates can make the structure functions constant or linear (homogeneous or non-homogeneously linear).
We wish to investigate the following conjecture of Michael Shapiro.

\begin{conj}
    If $f_1, \cdots, f_n \in \Rat$ are rational functions such that
    $\{f_i, f_j\} = \sum\limits_{k=1}^n c^k_{ij} f_k + d_{ij}$
    with $c^k_{ij}, d_{ij} \in \K$ for $1 \leq i,j,k \leq n$, then $\{f_i, f_j\} = 0$ for $1 \leq i,j \leq n$.
    \label{conj:linear}
\end{conj}

We prove this conjecture in Theorem~\ref{thm:main}.
Note that the conjecture implies that for any log-canonical Poisson structure on affine space, the answer to Question~\ref{q:degree} is ``no.''
That is, there is no system of coordinates whose structure
functions are polynomials of degree less than two.

\section{Nonexistence of Constant Bracket}
We have discovered that some of the results of this section have already appeared in~\cite{GL11}.
However, we include this section for completeness.
The results in this section will be built upon to prove our main theorem.
Given some $n \times n$ skew-symmetric matrix $\Omega = (\omega_{ij})$,  recall that $\Rat = \K(x_1, \cdots, x_n)$ 
is the algebra of rational functions in $n$ variables with the Poisson bracket given by
$$\{x_i, x_j\} = \omega_{ij} x_i x_j$$
for $1 \leq i,j \leq n$.
For $I = (i_1,\dots,i_n) \in \Z^n$, the corresponding Laurent monomial is written $\x^I = x_1^{i_1} \cdots x_n^{i_n}$. 
For $I = (i_1,\dots,i_n)$ and $J = (j_1,\dots,j_n)$ in $\Z^n$, let $A^I_J$ be the $2$-by-$n$ matrix whose rows are $I$ and $J$.
Let $\Delta_{ij}(A^I_J)$ be the $2$-by-$2$ minor of $A^I_J$ with columns indexed by $i$ and $j$.
Also define $M^I_J$ to be the weighted sum of the $\Delta_{ij}(A^I_J)$ given by the following formula
\[ M^I_J := \sum_{k < \ell} \omega_{k \ell} \Delta_{k \ell}(A^I_J) = \sum_{k < \ell} \omega_{k\ell} \left| \begin{array}{cc} i_k & i_\ell \\ j_k & j_\ell \end{array} \right| 
  = \sum_{k < \ell} \omega_{k\ell} (i_k j_\ell - i_\ell j_k). \]
Note that if $e_1,\dots,e_n$ is a basis for $\Bbb{Z}^n$, with $e^1,\dots,e^n$ the dual basis, we can define the two-form
\[ \omega = \sum_{k < \ell} \omega_{k\ell} \, e^k \wedge e^\ell, \] 
and then $M^I_J = \omega(I,J)$.
In particular, the expression $M^I_J$ is $\Bbb{Z}$-bilinear and skew-symmetric with respect to $I$ and $J$.
We now compute an explicit formula for the bracket of two Laurent polynomials.

\begin {lem}
    If $I, J \in \Z^n$, then
    $$ \{\x^I,\x^J\} = M^I_J \, \x^{I+J}.$$
    \label{lem:monomial}
\end {lem}
\begin{proof}
    Let $I = (i_1, \dots, i_n)$ and $J = (j_1, \dots, j_n)$.
    For $1 \leq k \leq n$ let $I_k = (i_1, \dots, i_{k-1}, 0, i_{k+1}, \dots, i_n)$ and $J_k =(j_1, \dots, j_{k-1}, 0, j_{k+1}, \dots, j_n)$.
    By using Equation \ref{eq:bracket-formula}, we find
    \begin {align*} 
        \{\x^I, \x^J\} &= \sum_{1 \leq k, \ell \leq n} \frac{\partial \x^I}{\partial x_k} \frac{\partial \x^J}{\partial x_\ell} \{x_k,x_\ell\} \\
                     &= \sum_{1 \leq k, \ell \leq n} i_k \, j_\ell \, \x^{I_k+J_{\ell}} \{x_k, x_{\ell}\} \\
                     &= \sum_{1 \leq k, \ell \leq n} i_k \, j_\ell \, \x^{I_k+J_{\ell}} \omega_{k \ell} x_k \, x_{\ell} \\
                     &= \sum_{1 \leq k, \ell \leq n} \omega_{k \ell} i_k \, j_\ell \, \x^{I+J} \\
                     &= \sum_{1 \leq k < \ell \leq n} \omega_{k \ell}(i_k j_\ell - i_\ell j_k) \x^{I+J} \\
                     &= M^I_J \x^{I+J}
    \end {align*}
\end{proof}

In order to prove our first theorem we want to work with iterated Laurent series.
We will give a brief overview of the theory of iterated Laurent series which will be needed for our purpose.
For a more in depth treatment of iterated Laurent series we refer the reader to~\cite{Xin04}.
For us a \emph{formal Laurent series} in variables $x_1, \dots, x_n$ over $\K$ will mean any formal sum
$$f = \sum_{I \in \Z^n} \alpha_I \x^I$$
with $\alpha_I \in \K$ for all $I \in \Z^n$.
For any $I \in \Z^n$ let $[\x^I]f$ denote the coefficient of $\x^I$ in $f$.
In particular, $[\mathbf{1}]f$ denotes the constant term of $f$.
Also, we let $\supp(f)$ denote the set $I \in \Z^n$ such that $[\x^I]f \ne 0$.
The set of formal Laurent series is a $\K$-vector space, but it is not a $\K$-algebra as we cannot multiply any two formal Laurent series in general.

However, certain subsets of the set of formal Laurent series form a $\K$-algebra.
Define $\K\langle\langle x \rangle\rangle := \K((x))$ to be the field of Laurent series in a single variable.
That is, $\Bbb{K}((x))$ consists of formal Laurent series $\sum_{i \in \Bbb{Z}} a_i x^i$ containing only finitely many negative exponents.
Now define $$\K\langle \langle x_1, \cdots, x_{i+1}\rangle\rangle := \K\langle\langle x_1, \cdots, x_i \rangle\rangle((x_{i+1}))$$ iteratively.
We then let $\Lau = \K \langle\langle x_1, \cdots, x_n\rangle \rangle$ be the field of \emph{iterated Laurent series} in $n$ variables.
We have the following immediate corollary of Lemma~\ref{lem:monomial} which holds for Laurent polynomials.
In the remainder of this section we will show that this corollary can be extended to hold for any iterated Laurent series.

\begin {cor}
    Let $f,g \in \Bbb{K}[x_1^{\pm},\dots,x_n^{\pm}]$ be Laurent polynomials, with $\mathcal{I} = \supp(f)$ and $\mathcal{J} = \supp(g)$.
    Then their bracket is given by
    \[ \{f,g\} = \sum_{(I,J) \in \mathcal{I} \times \mathcal{J}} \alpha_I \beta_J M^I_J \x^{I+J}. \]
    \label{cor:polynomial}
\end {cor}

\begin {rem}
    Note that we have the inclusion $\K[x_1, \dots, x_n] \hookrightarrow \Lau$.
    Since $\Lau$ is a field and $\Rat$ is the field of fractions of $\K[x_1, \dots, x_n]$, 
    we also have the inclusion $\Rat \hookrightarrow \Lau$.
    Hence, $\Rat$ is a $\Bbb{K}$-subalgebra of $\Lau$.
\end {rem}

\bigskip

\begin {rem}
    Notice the order in which we adjoin our variables is relevant. For instance, consider the rational function $\frac{1}{x+y}$.
    As an element of $\K \langle\langle x,y \rangle\rangle$, it can be written as
    \[ \frac{1}{x+y} = \sum_{n \geq 0} (-1)^n x^{-(n+1)}y^n \]
    However, since there is no lower bound on the powers of $x$, this does not give an element of $\K \langle\langle y,x \rangle\rangle$.
    Instead, to represent it as an element in the latter field, we must write
    \[ \frac{1}{x+y} = \sum_{n \geq 0} (-1)^n x^ny^{-(n+1)}. \]
\end {rem}

\bigskip

\begin {rem}
    Any iterated Laurent series $f \in \Lau$ can be expressed as a formal Laurent series.
    That is, we can write
    $$f = \sum_{I \in \supp(f)} \alpha_I \x^I.$$
    Given $f, g \in \Lau$ where
    \begin{align*}
    f &= \sum_{I \in \supp(f)} \alpha_I \x^I & g &= \sum_{J \in \supp(g)} \beta_J \x^J
    \end{align*}
    their product is
    $$fg = \sum_{(I,J) \in \supp(f) \times \supp(g)} \alpha_I \beta_J \x^{I+J}.$$
    This product $fg$ is also an iterated Laurent series, since $\Lau$ is a field.
    In particular this means $fg$ is a formal Laurent series with the property that for any $K \in \Z^n$, the set
    $$\{(I,J) \in \supp(f) \times \supp(g) : I+J=K\}$$
    is finite. In fact, we have the following result, which will be useful later\footnote{We have chosen to use the iterated Laurent construction, and hence must show the well-ordered support property in Proposition~\ref{prop:well-ordered}. Alternatively, we could have started from the well-ordered support property and shown that we obtain a ring structure. A formal series with well-order support are sometimes called a \emph{Hahn series} or a \emph{Mal'cev-Neumann series} and exponents can be taken from any ordered abelian group.}:
\end {rem}

\bigskip

\begin {prp}[{\cite[Proposition 2-2.1]{Xin04}}]
    Let $f$ be a formal Laurent series. Then $f \in \Lau$ if and only if
    $\supp(f)$ is well-ordered with respect to the reverse lexicographic ordering.
    \label{prop:well-ordered}
\end {prp}

\bigskip

\begin{lem}
    The Poisson bracket on $\Rat$ extends uniquely to a Poisson bracket on $\Lau$.
\label{lem:sub}
\end{lem}
\begin{proof}
    Note that by bilinearity, any Poisson bracket on $\Lau$ is determined by the brackets of all Laurent monomials.
    Thus by Lemma~\ref{lem:monomial}, any bracket extending the one on $\Rat$ must be given by the same formula on monomials.
    We claim that the same formula in Corollary \ref{cor:polynomial} gives the bracket on $\Lau$.
    It suffices to show that for $f,g \in \Lau$ that $\{f,g\} \in \Lau$. That is we must
    show that given $f,g \in \Lau$, the formula from Corollary \ref{cor:polynomial} yields an element of $\Lau$.

    Let $f, g \in \Lau$, and again use the notation $\mathcal{I} = \supp(f)$ and $\mathcal{J} = \supp(g)$.
    Note that since $fg \in \Lau$, then $\supp(fg)$ is well-ordered by Proposition \ref{prop:well-ordered}.
    The formula from Corollary \ref{cor:polynomial} also implies that $\supp(\{f,g\}) \subseteq \supp(f) + \supp(g)$,
    where ``$+$'' is used to denote Minkowski addition:
    \[ \supp(f) + \supp(g) = \{I+J ~|~ I \in \supp(f), \, J \in \supp(g) \}. \]
    Being a subset of a well-ordered set, we see that $\supp(\{f,g\})$ is itself well-ordered. Hence, $\{f,g\} \in \Lau$ by
    Proposition~\ref{prop:well-ordered}.

\end{proof}

The remaining results in this section are restatements of the indicated results from~\cite{GL11}. The next theorem is a simple but powerful observation
which is an essential ingredient to our proof of Conjecture~\ref{conj:linear}.

\begin{thm}[{\cite[Proposition 5.2 (a)]{GL11}}]
    If $f, g \in \Lau$, then $[\mathbf{1}]\{f, g\} = 0$.
    \label{thm:constant}
\end{thm}
\begin{proof}
    As usual, let $\I$ and $\J$ be the supports of $f$ and $g$, and let
    \begin{align*}
    f &= \sum_{I \in \I} \alpha_I \x^I & g &= \sum_{J \in \J} \beta_J \x^J   
    \end{align*}
    be expressions for $f$ and $g$ as formal Laurent series.
    Computing using Corollary~\ref{cor:polynomial} we see that
    $$\{f,g\} = \sum_{(I,J) \in \I \times \J} \alpha_I \beta_J M^I_J \x^{I+J}$$
    and so
    $$[\mathbf{1}]\{f,g\} = \sum_{\substack{(I,J) \in \I \times \I \\ I+J = 0}} \alpha_I \beta_J M^I_J.$$
    However, if $I + J = 0$, then $J = -I$ and $M^I_J = 0$.
    Here we have used that $M^I_J$ is skew-symmetric with respect to $I$ and $J$.
\end{proof}

\bigskip

\begin{cor}[{\cite[Corollary 5.3]{GL11}}]
    If $f_1, \cdots, f_n \in \Rat$ are rational functions such that
    $\{f_i, f_j\} = c_{ij}$
    with $c_{ij} \in \K$ for $1 \leq i,j \leq n$, then $c_{ij} = 0$ for $1 \leq i,j \leq n$.
    \label{cor:rational_constant}
\end{cor}
\begin{proof}
    By Lemma~\ref{lem:sub}, $\Rat$ is a Poisson subalgebra of $\Lau$.
    The corollary then follows immediately from Theorem~\ref{thm:constant}.
\end{proof}

\section{Nonexistence of Linear Bracket}

As in the previous section, we consider the Poisson algebra of rational functions $\Rat$, in $n$ variables, with bracket given by
\[ \{x_i,x_j\} = \omega_{ij} x_i x_j \]
for some skew-symmetric matrix $\Omega = (\omega_{ij})$ with coefficients in $\K$. It is the goal of this section to prove the aforementioned Conjecture~\ref{conj:linear},
which states that there is no rational change of coordinates making the bracket linear. That is, if there are rational functions $f_1,\dots,f_n$
such that $\{f_i,f_j\} = \sum_{k=1}^n c^k_{ij} f_k + d_{ij}$ for constants $c^k_{ij}, d_{ij} \in \Bbb{K}$, then in fact all the coefficients $c^k_{ij}$ and $d_{ij}$
must be zero.

\bigskip

We now prove a lemma which will be used later.

\begin {lem}
    There do not exist linearly independent $f,g \in \Rat$ such that $\{f,g\} = af+bg$ for $a,b \in \K$ with $a$ and $b$ not both zero.
    \label{lem:linear}
\end {lem}
\begin {proof}
    Assume there do exist linearly independent rational functions $f$ and $g$ so that $\{f,g\} = af+bg$ for some $a, b \in \K$. Then the linear span of $f$ and $g$ is a two-dimensional
    Lie subalgebra of $\Rat$. Up to isomorphism, there is a unique two-dimensional non-abelian Lie algebra, with the bracket given
    by $\{f,g\} = f$. Explicitly, the isomorphism is given by $f \mapsto f + \frac{b}{a}g$, $g \mapsto \frac{1}{a} g$ (assuming $a \neq 0$).
    So, we may assume without loss of generality that $a=1$ and $b=0$, thus $\{f,g\} = f$. But then we have that $\frac{1}{f}\{f,g\} = 1$.
    Note that since $\mathrm{ad}_f = \{f,\cdot\}$ is a derivation, we have $\frac{1}{f}\{f,g\} = \left\{f, \frac{g}{f} \right\}$. This in turn
    implies that $\{f,\frac{g}{f}\} = 1$. But this directly contradicts Corollary~\ref{cor:rational_constant}, which says that the bracket of any two rational functions
    cannot be a nonzero constant.
\end {proof}

\bigskip

A useful consequence of this lemma is that the adjoint maps $\ad_f$ can have no non-zero eigenvalues.

\begin{cor}
    If $f,g \in \Rat$ with $g \neq 0$ and $\{f,g\} = \lambda g$ for some $\lambda \in \K$, then $\lambda = 0$.
    \label{cor:eigen}
\end{cor}

\bigskip

The next result says that the adjoint maps $\ad_f$ cannot be nonzero and nilpotent.

\begin{lem}
    If $f,g \in \Rat$ and $\{f,g\} \neq 0$, then $\{f, \{f,g\}\} \neq 0$.
    \label{lem:nil}
\end{lem}

\begin{proof}
Take $f,g \in \Rat$ and assume that $\{f,g\} \neq 0$ but $\{f, \{f,g\}\} = 0$.
Then we know that $\left\{f, \frac{1}{\{f,g\}} \right\} = 0$.
Computing, we see that
$$\left\{f, \frac{g}{\{f,g\}}\right\} = g\left\{f, \frac{1}{\{f,g\}}\right\} + \frac{1}{\{f,g\}} \{f,g\} = 1$$
which is a contradiction to Corollary~\ref{cor:rational_constant}.
\end{proof}

\bigskip

We are now ready to prove the main result.
\begin{thm}
    If $f_1, \cdots, f_n \in \Rat$ are rational functions such that
    $$\{f_i, f_j\} = \sum\limits_{k=1}^n c^k_{ij} f_k + d_{ij}$$
    with $c^k_{ij}, d_{ij} \in \K$ for $1 \leq i,j,k \leq n$, then $\{f_i,f_j\} = 0$ for  $1 \leq i,j\leq n$.
\label{thm:main}
\end{thm}

\begin {proof}
    Assume first that $\K = \overline{\K}$ is algebraically closed. 
    Let $f_1,\dots,f_n \in \Rat$ be rational functions such that $\{f_i,f_j\} = \sum_{k=1}^n c^k_{ij} f_k + d_{ij}$ for some $c^k_{ij}, d_{ij} \in \Bbb{K}$.
    This means that $1, f_1,\dots,f_n$ generate a finite dimensional Lie algebra inside $\Rat$.
    Let $F \leq \Rat$ denote this finite dimensional Lie algebra generated by $1,f_1, \dots, f_n$.
    For any $f \in F$ we have the linear map $\ad_f: F \to F$, and by Corollary \ref{cor:eigen} all eigenvalues of $\ad_f$ are zero.
    It follows, since $\K$ is algebraically closed, that $\ad_f$ is nilpotent.
    However, Lemma~\ref{lem:nil} implies that if $\ad_f$ is nilpotent we must have $\ad_f = 0$.
    The theorem then follows.

    In the case that $\K$ is \emph{not} algebraically closed, consider $\overline{\Rat} := \overline{\K} \otimes_{\K} \Rat = \overline{\K}(x_1,\dots,x_n)$.
    The relations $\{f_i,f_j\} = \sum_{k=1}^n c^k_{ij} f_k + d_{ij}$ still hold.
    Thus $1,f_1, \dots, f_n$  will generate some finite dimensional Lie algebra inside $\overline{\Rat}$, and we can complete the argument just as in the algebraically closed case.
\end {proof}

\begin{rem}
Given a Poisson algebra $P$, the \emph{quadratic Poisson Gel'fand-Kirillov problem} is to determine if the field of fractions of $P$ is isomorphic to $\Rat$ for some $\Omega$.
This problem was first defined in~\cite{GL11} and further studied in~\cite{LL16}.
In this section, we have shown a number of properties of the Poisson algebra $\Rat$.
Hence, any Poisson algebra isomorphic to $\Rat$ must also have these properties, and the results in this section 
can be viewed as necessary conditions for a Poisson algebra to be a solution to the quadratic Poisson Gel’fand-Kirillov problem.
\end{rem}

\section{Generalizations}

The results of the previous section are not specific to only the Poisson algebra $\Rat$.
Let $P$ be a Poisson algebra $P$ with the following two properties:
\begin{itemize} 
    \item $P$ is a field.
    \item For any $a, b \in P$ we have $\{a,b\} = 0$ whenever $\{a,b\} \in \K$.
\end{itemize}
Call such an algebra a \emph{nonconstant Poisson field} (since there are no elements with $\{f,g\} = 1$).
Then versions of the results in the previous section hold for $P$ since the proofs only use the conditions above.
In particular we will have a version of Theorem~\ref{thm:main} which says that $P$ can have no finite dimensional non-abelian Lie subalgebra. 
Before proving this theorem, let us collect some of the essential parts of the proofs from the previous section into a useful general lemma:

\begin {lem}
    Let $P$ be a Poisson $\Bbb{K}$-algebra which is a field. Then the following are equivalent:
    \\ \\
    \begin {tabular}{cp{5in}}
        $(a)$ & There exist $f,g \in P$ such that $\{f,g\} = 1$. \\ [1.5ex]
        $(b)$ & There exist $f,g \in P$ such that $\{f,g\} = g$. \\ [1.5ex]
        $(c)$ & There exist $f,g \in P$ with $\{f,g\} \neq 0$ but $\{f,\{f,g\}\} = 0$.
    \end {tabular}
\label{lem:gen}
\end {lem}
\begin {proof} \ \\ \\ \relax
    $(b) \Rightarrow (a)$: Follows from proof which is identical to the proof of Lemma \ref{lem:linear}.

    \bigskip

    \noindent
    $(c) \Rightarrow (a)$: Follows from proof which is identical to the proof of Lemma \ref{lem:nil}

    \bigskip

    \noindent
    $(a) \Rightarrow (c)$: If $\{f,g\} = 1$, then $\{f,\{f,g\}\} = \{f,1\} = 0$.

    \bigskip

    \noindent
    $(a) \Rightarrow (b)$: Suppose that $\{f,g\} = 1$, and define $x = fg$ and $y = g$. Then 
    \[ \{x,y\} = \{fg,g\} = \{f,g\}g = g = y. \]
\end {proof}

\bigskip

Analogous to the above definition, define a \emph{nonlinear Poisson field} as a Poisson field $P$ which
has no finite-dimensional nonabelian Lie subalgebras. This means there are no finite collections $f_1,\dots,f_k \in P$
and constants $c^\ell_{ij}$ such that $\{f_i,f_j\} = \sum_\ell c^\ell_{ij} f_\ell$.
The next result says that for a Poisson field, being nonconstant is a sufficient condition to being nonlinear.

\begin {thm}
    Any nonconstant Poisson field is also nonlinear.
\end {thm}
\begin {proof}
We assume that we are working over an algebraically closed field, if not we can modify just as in the proof of Theorem~\ref{thm:main}.
    Suppose that there exist some $f_1,\dots,f_k \in P$ for some $k > 1$ and constants $c^\ell_{ij}$ so that
    \[ \{f_i,f_j\} = \sum_\ell c^\ell_{ij} f_\ell \]
    Then $f_1,\dots,f_k$ generate a finite-dimensional Lie subalgebra $F \leq P$. Each map $\mathrm{ad}_{f_i}$ is an endomorphism of $F$.
    Note that $\mathrm{ad}_{f_i}$ cannot have any nonzero eigenvalues. If it did, there would be
    some $g \in F$ and $\lambda \neq 0$ so that $\{f_i,g\} = \lambda g$. Then for $\tilde{f}_i = \frac{1}{\lambda} f_i$, we have $\{\tilde{f}_i,g\} = g$.
    By the previous theorem, there must also exist some $u,v \in P$ so that $\{u,v\} = 1$. But this contradicts the assumption on $P$.
    So in fact $\mathrm{ad}_{f_i}$ can have only zero eigenvalues, and hence must be nilpotent. Again, by the previous theorem, if
    $\mathrm{ad}_{f_i}$ is nonzero and nilpotent, then there would exist $u,v \in P$ with $\{u,v\} = 1$. So it must be that $\mathrm{ad}_{f_i} = 0$,
    and thus $F$ is an abelian Lie algebra.
\end {proof}

In the spirit of Question~\ref{q:degree}, let us call a system of coordinates \emph{(homogeneous) algebraically reduced} 
if all structure functions are (homogeneous) polynomials of a given degree, and there does not exist any rational change of coordinates 
making the structure functions (homogeneous) polynomials of a smaller degree.
In Theorem~\ref{thm:main} we provided an answer to Question~\ref{q:degree} for any log-canonical system of coordinates and showed that 
they are algebraically reduced.
It is natural to look for other (homogeneous) algebraically reduced coordinate systems.

\bigskip 

Let us now consider systems of coordinates for which all structure functions are monomials.
In dimension $2$ with coordinates $(x,y)$ so that $\{x,y\} = x^ay^b$ we have seen that such a monomial system of coordinates is algebraically reduced if and only if $(a,b) = (0,0)$ or $(a,b) = (1,1)$.
In dimension $3$ with coordinates $(x,y,z)$ and bracket relations
\begin{align*}
\{x,y\} &= Ax^{a_1}y^{a_2}z^{a_3}\\
\{x,z\} &= Bx^{b_1}y^{b_2}z^{b_3}\\
\{y,z\} &= Cx^{c_1}y^{c_2}z^{c_3}
\end{align*}
we can extend by skew-symmetry, but must also ensure the Jacobi identity holds.
Computing we obtain
\begin{align*}
\{x,\{y,z\}\} + \{y,\{z,x\}\} + \{z,\{x,y\}\} &= (b_1 - a_1)AB x^{a_1 + b_1 -1} y^{a_2 + b_2} z^{a_3 + b_3}\\
&\quad\quad + (c_2 - a_2)AC x^{a_1 + c_1} y^{a_2 + c_2 - 1} z^{a_3 + c_3}\\
&\quad\quad + (c_3 - b_3)BC x^{b_1 + c_1} y^{b_2 + c_2} z^{b_3 + c_3 - 1}.
\end{align*}
If $a_1 = b_1, a_2 = c_2,$ and $b_3 = c_3$ the Jacobi identity will hold.
In that case the bracket relations are
\begin{align*}
\{x,y\} &= A(x^{a_1}y^{a_2}z^{b_3})z^{a_3 - b_3}\\
\{x,z\} &= B(x^{a_1}y^{a_2}z^{b_3})y^{b_2 - a_2}\\
\{y,z\} &= C(x^{a_1}y^{a_2}z^{b_3})x^{c_1 - a_1}
\end{align*}

Consider the simplest example of the case above, where $(a_1,a_2,b_3) = (0,0,0)$. One such example is the bracket on $\K(x,y,z)$ given by
\begin{align*}
\{x,y\} &= a z^2 &\{x,z\} &= b y^2 & \{y,z\} &= c x^2
\end{align*}
for some $a,b,c \in \K^*$.
This bracket gives a candidate for another homogeneous quadratic algebraically reduced system of coordinates which differs from the log-canonical case.
By the above discussion, it would suffice to show that this bracket makes $\K(x,y,z)$ a nonconstant Poisson field.
However, unlike the log-canonical case, this bracket can produce non-zero constant terms, as exhibited by the following examples:
$$\left\{x, \frac{y}{z^2}\right\} = a - 2b\left(\frac{y}{z}\right)^3$$
$$\left\{\frac{x}{z}, \frac{y}{z}\right\} = a - b\left(\frac{y}{z}\right)^3 + c\left(\frac{x}{z}\right)^3$$
As such, the arguments used in the present paper do not apply, since everything followed from Theorem \ref{thm:constant},
which said that the constant term of $\{f,g\}$ (viewed as a Laurent series) is always zero. However, it is possible that this
bracket makes $\K(x,y,z)$ a nonconstant Poisson field, despite the fact that Theorem \ref{thm:constant} does not hold.

\bigskip

It seems to be an interesting problem to find other algebraically reduced brackets on $\K(x_1,\dots,x_n)$,
and to find necessary and/or sufficient conditions on the structure functions.

\bibliographystyle{alpha}
\bibliography{PoissonBib}
\end{document}